\documentclass[12pt]{article} \textwidth=6in \oddsidemargin=0in
\usepackage{amssymb} 
\begin{document}\def\ov{\over} 
\newcommand{\C}[1]{{\cal C}_{#1}} \def\inv{^{-1}}\def\be{\begin{equation}} \def\ep{\varepsilon} 
\def\ee{\end{equation}}\def\x{\xi}\def\({\left(} \def\){\right)} \def\iy{\infty}  \def\ld{\ldots} \def\Z{\mathbb Z} \def\cd{\cdots}
\def\P{\mathbb P}\def\bs{\backslash} \def\t{\tau} 
\def\r{\rho}  \def\ph{\varphi}  \def\e{\eta} 
\newcommand{\br}[2]{\left[{#1\atop #2}\right]_\t}  \def\la{\lambda} \def\sp{\vspace{2ex}}  \def\l{\ell} \def\s{\sigma} \def\tl{\tilde} 
\def\la{\lambda} \newcommand{\brr}[2]{\left[{#1\atop #2}\right]} 
\def\DI{\Delta I} \def\mb{\boldmath}

\hfill   January 26, 2010

\begin{center}{\bf \large Formulas for ASEP with Two-Sided\\ 
\vspace{.5ex}Bernoulli Initial Condition}\end{center}

\begin{center}{\large\bf Craig A.~Tracy}\\
{\it Department of Mathematics \\
University of California\\
Davis, CA 95616, USA\\
email: tracy@math.ucdavis.edu}\end{center}

\begin{center}{\large \bf Harold Widom}\\
{\it Department of Mathematics\\
University of California\\
Santa Cruz, CA 95064, USA\\
email: widom@ucsc.edu}\end{center}
 \begin{abstract}
For the asymmetric simple exclusion process $\e_t$ on the integer lattice  with two-sided Bernoulli initial condition, we derive
exact formulas for the following quantities: (1)  $\P(\e_t(x)=1)$, the probability that site $x$ is occupied at time~$t$; (2) the correlation
function $\P(\e_t(x)=1, \e_0(0)=1)$; (3) the distribution function for $Q_t$, the total flux across $0$ at time~$t$, and its exponential generating function.
 \end{abstract}
\begin{center}{\bf I. Introduction}\end{center}

Since its introduction by Spitzer \cite{S}, the asymmetric simple exclusion process (ASEP) has attracted considerable attention in both the mathematics and physics literature due to the fact it is one of the simplest lattice models describing transport far from equilibrium \cite{D, L1, L2}.  Recall that ASEP
on the integer lattice $\Z$ is a continuous time Markov process $\e_t$, where $\e_t(x)=1$ if $x\in\Z$
is occupied at time $t$, and $\e_t(x)=0$ if $x$ is vacant at time $t$.   Particles move on $\Z$ according to two rules: (1) A particle at $x$ waits an
exponential time with parameter one (the particle's  ``alarm clock''), and then chooses $y$ with probability $p(x,y)$; (2) If $y$ is vacant at that time it moves to $y$, while if $y$ is occupied it remains at $x$ and restarts its clock.  The process is called ``simple''
since the jumps are restricted to nearest neighbors: $p(x,x+1)=p$ and $p(x,x-1)=q=1-p$.  
All  clocks are independent of each other and one need not worry that two clocks ring simultaneously since this event has probability zero.  A rigorous construction of this process can be found in Liggett \cite{L1}.  We note that the process is called the totally asymmetric simple exclusion process (TASEP) if jumps occur only to the left ($q=1$) or only to the right ($p=1$); and if $p=q=1/2$, the process is called the symmetric simple exclusion process (SSEP) \cite{L1, L2}.

In previous work \cite{TW3} the authors considered ASEP on the integer lattice $\Z$ with  \textit{step Bernoulli initial condition}: initially a site in $\Z^+$ (the positive integers) is occupied with probability $\r$ while sites in $\Z^-$ (the nonpositive integers) are unoccupied. 
Formulas were obtained for $\P(x_m(t)\le x)$, the probability that at time $t$ the $m$th particle from the left is at a site $\le x$, with consequent asymptotic results for this probability.  In this paper we consider the more general situation where initially a site in $\Z^+$ is occupied with probability $\r_+$ while a site in $\Z^-$ is occupied with probability $\r_-$.  This initial condition is called the  \textit{two-sided Bernoulli initial condition}.  Now it makes no sense to speak of the $m$th particle from the left. Instead we find formulas for $\P(\e_t(x)=1)$, the probability that site $x$ is occupied at time $t$, and the correlation function $\P(\e_t(x)=1,\,\e_0(0)=1)$. We also derive formulas for the generating function $\langle e^{\la Q_t}\rangle$, where $Q_t$ is the total flux across $0$ at time $t$.  

For the special case of TASEP with two-sided Bernoulli initial condition much is already known including various limit theorems \cite{BAC, PS}.  (Also see \cite{FS}
for the stationary case $\r_{-}=\r_{+}$.) A simplifying feature of TASEP is that it is a determinantal process \cite{J} thus permitting the application of random matrix theory techniques.  In \cite{TW3} it was shown that the limit theorems  in \cite{BAC} for TASEP extend to ASEP for the special case of
step Bernoulli initial condition ($\r_{-}=0$, $\r_{+}>0$); and thus, these results \cite{TW3} establish a KPZ universality theorem for ASEP.
Whether the TASEP limit theorems of \cite{BAC} for two-sided Bernoulli initial condition can be extended to ASEP from the
formulas derived in this paper remains to be seen.  

We shall assume throughout that $p,\,q\ne0$. As in \cite{TW3} we begin with a formula derived in \cite{TW1} for $\P_Y(x_m(t)\le x)$, the probablity distribution for the $m$th particle from the left given an initial finite configuration $Y$; then we average over all $Y$ in a finite subset of $\Z$; and then we let the subset of $\Z$ become unbounded. The formula for $\P_Y(x_m(t)\le x)$ involved multiple integrals over large contours, and this allowed $Y$ to become unbounded on the right in \cite{TW3}. Now, if $Y$ is to become unbounded on left and right, we must sum over both large and small contours, so the formula for $\P_Y(x_m(t)\le x)$ must first be recast. This is done in Sec.~II. In Sec. III we deduce the formula for $\P_Y(\e_t(x)=1)$, the probability that site $x$ is occupied at time $t$ given the initial configuration $Y$. In Sec. IV we average over $Y$ and take the limit to obtain the formula for $\P(\e_t(x)=1)$ for two-sided Bernoulli initial configuration on $\Z$. The formula for the correlation function requires only a small change. Finally, in Sec.~V we deduce formulas for probability distribution of the total flux across the origin at time $t$ and for its exponential generating function. For the special case of SSEP we connect this last result with a recent paper of Derrida and Gerschenfeld~\cite{DG}. 

\begin{center}{\bf II. Finite Configuration, Large and Small Contours}\end{center}

We begin with Theorem 5.2 of \cite{TW1} which gives a formula for $\P_Y(x_m(t)\le x)$, the probability in ASEP when the initial configuration is a (deterministic) finite set $Y$. To state it we introduce some notation, which here will be slightly different. 

First, we set $\t=p/q$ and recall that the $\t$-binomial coefficient $\br{N}{n}$ is defined as
\[{(1-\t^N)\,(1-\t^{N-1})\cdots (1-\t^{N-n+1})\ov (1-\t)\,(1-\t^2)\cdots (1-\t^n)}\]
when $n$ is a positive integer, 1 when $n=0$ (empty products are always defined to be~1), and 0 when $n$ is a negative integer.  

We define
\[\ep(\x)=p\,\x\inv+q\,\x-1,\ \ \ \ f(\x_i,\,\x_j)={\x_j-\x_i\ov p+q\x_i\x_j-\x_i},\] 
and then
\[I(x,k,\x)=I(x,k,\x_1,\ld,\x_k)=\prod_{i<j}f(\x_i,\,\x_j)\;\prod_i{\x_i^{x}\,e^{\ep(\x_i)t}\ov 1-\x_i},\]
\[\DI(x,k,\x)=I(x,k,\x)-I(x-1,k,\x)=(1-\prod_i\x_i\inv)\;I(x,k,\x).\]
All indices in the products run from 1 to $k$. Notice that $I$ and $\DI$ depend on $t$, although it is not displayed in the notations.

Finally, given two sets of integers $U$ and $V$ we define
\[\s(U,\,V)=\#\{(u,\,v): u\in U,\ v\in V,\ {\rm and}\ u\ge v.\}.\]
(If $V=[1,\,N]$\footnote{This is short for the set $\{1,\ld,N\}$. We use this notation consistently.} and $U\subset [1,\,N]$ then $\s(U,\,V)$ equals the sum of the elements of $U$. Hence the notation.)

Theorem 5.2 of \cite{TW1} states that when $q\ne0$,
\be\P_Y(x_m(t)=x)=\sum_{k=1}^{|Y|}c_{m,k}\,\sum_{{S\subset Y\atop |S|=k}}\,\t^{\s(S,\,Y)}\,\int_{\C{R}}\cd\int_{\C{R}} \DI(x,k,\x)\,\prod_i\x_i^{-s_i}\,d\x_1\cd d\x_k,\label{PY1}\ee
where
\[S=\{s_1,\ld,s_k\},\]
\be c_{m,k}=(-1)^{m}\,q^{k(k-1)/2}\,\t^{m(m-1)/2-km}\,
\br{k-1}{m-1},\label{c}\ee
and $\C{R}$ is the circle with center 0 and radius $R$, which is assumed so large that the denominators $p+q\x_i\x_j-\x_i$  are nonzero on and outside the contours. The inner sum is taken over all subsets $S$ of $Y$ with $|S|=k$. Observe that $c_{m,k}=0$ when $k<m$.

If we sum (\ref{PY1}) on $x$ from $-\iy$ to $x$, as we may do since $R>1$, we obtain
\be\P_Y(x_m(t)\le x)=\sum_{k\ge1}c_{m,k}\,\sum_{{S\subset Y\atop |S|=k}}\,\t^{\s(S,\,Y)}\,\int_{\C{R}}\cd\int_{\C{R}} I(x,k,\x)\,\prod_i\x_i^{-s_i}\,d\x_1\cd d\x_k.\label{PY2}\ee 
Clearly (\ref{PY2}) implies (\ref{PY1}) immediately .

Because of the appearance of the factors $\prod_i\x_i^{-s_i}$ in the integrands and because $R$ is large we expect to be able to allow $Y$ to be unbounded on the right but not on the left. So we want integrals over both large and small contours. Lemma 5.1 of \cite{TW1} gives a formula which expresses an integral over a large contour as a sum of integrals over small contours, and we shall be able to apply it here. 

We shall assume throughout that $Y=Y_-\cup Y_+$ where $Y_-$ is to the left of $Y_+$. Then each $S\subset Y$ is a union $S_-\cup S_+$ with $S_\pm\subset Y_\pm$. Set $k_\pm=|S_\pm|$, so $k=k_-+k_+$. Instead of using indices $i=1,\ld,k$ for each $S$ with $|S|=k$, we shall use $k_-$ negative indices $-k_-,\ld,-1$ and $k_+$ positive indices $1,\ld,k_+$, with $s_i\in S_-$ if $i<0$ and $s_i\in S_+$ if $i>0$. In this notation (\ref{PY2}) becomes
\[\P_Y(x_m(t)\le x)=\sum_{k_\pm\ge0,\ k\ge1}c_{m,k}\,\sum_{{S_\pm\subset Y_\pm\atop |S_\pm|=k_\pm}}\,\t^{\s(S_-,\,Y_-)+\s(S_+,\,Y)}\]
\be\times\int_{\C{R}}\cd\int_{\C{R}}I(x,k,\x)\,\prod_i\x_i^{-s_i}\,\prod_i d\x_i\label{PY3},\ee
where $k=k_-+k_+$. For the exponent of $\t$ we used the bilinearity of $\s$ and the obvious fact $\s(S_-,Y_+)=0$.

We rephrase Lemma 5.1 of \cite{TW1} to make it compatible with the present notation. Let $A$ be a finite set of indices and let
\[g=g(\x_i)_{i\in A}\] 
be a function that is analytic for all $\x_i\ne0$. Assume that for $i<k$
\[g\Big|_{\x_i\to p/(1-q\x_k)}=O(\x_k\inv)\]
as $\x_k\to\iy$, uniformly when all the $\x_\l$ with $\l\ne i,\,k$ are bounded and bounded away from zero. For $B\subset A$ denote by $g_B$ the function obtained from $g$ by setting all $\x_i$ with $i\not\in B$ equal to 1. (In particular $g=g_A$.) Define
\[I_B(\x)=\prod_{{i<j\atop i,j\in B}}f(\x_i,\x_j)\,{g_B(\x)\ov\prod_{i\in B}(1-\x_i)}.\]
Then when $p,\,q\ne0$,
\be\int_{\C{R}^{|A|}} I_A(\x)\,\prod_{i\in A}d\x_i=\sum_{B\subset A}
(-1)^{|A\bs B|}\t^{\s(B,\,A)-|A|\,|B|}\;{p^{|B|(|B|-1)/2}\ov q^{|A|(|A|-1)/2}}\,
\int_{\C{r}^{|B|}} I_B(\x)\,\prod_{i\in B}d\x_i,\label{gint}\ee
where $r$ is so small that all the zeros of the denominators lie outside $\C{r}$. When $B$ is empty the integral on the right side is interpreted as $g(1,\ld,1)$.

We shall apply this to the integral on the right side of (\ref{PY3}). 
In order to do this we change the set of indices in that integral from  $[-k_-,\,-1]\cup[1,\,k_+]$ to $S_-\cup S_+$ in the obvious way. In the application of the lemma we take $A$ to be the set $S_-$ of negative indices, and $g$ the integral over the remaining variables $\x_i$ with $i>0$ (i.e., $i\in S_+$):
\[g(\x)=\int_{\C{R}^{k_+}}\prod_{{i<j\atop j>0}}f(\x_i,\,\x_j)\;{\prod_i\x_i^{x-s_i}\,e^{\ep(\x_i)t}\ov \prod_{i>0}(1-\x_i)}\;\prod_{i>0}\,d\x_i.\]
Then the integral on the right side of (\ref{PY3}) equals the left side of (\ref{gint}).

Let us see that $g(\x)$ satisfies the required conditions. Since $R$ is arbitrarily large, $g(\x)$ is analytic for $\x_i\ne 0\ (i<0)$.  When $i<k<0$ and we set $\x_i=p/(1-q\x_k)$ the terms in the integral defining $g$ that involve $\x_k$ are 
\[\prod_{j>0}f\({p\ov1-q\x_k},\,\x_j\)
\,\cdot\,\prod_{j>0}f(\x_k,\,\x_j)\,\cdot\,\({p\ov1-q\x_k}\)^{x-s_i}\,\x_k^{x-s_k}\,e^{\,[p/\x_k+pq/(1-q\x_k)]t}.\]
As $\x_k\to\iy$ the product of the $f$s as well as the last factor are clearly bounded as $\x_k\to\iy$. The product of the remaining two factors is $O(\x_k^{s_i-s_k})\to0$ since $s_i<s_k$.

So the hypothesis of the lemma holds. In its application\footnote{In the resulting integrals over $\C{r}$ we must have $rR$ large, since to maintain the analyticity of $g(\x)$ for small $\x_i$ (with $i<0$) we must have all the $\x_j$ (with $j>0$) even larger.\label{rR}} we replace $|A\bs B|=k_--|B|$ of the $\x_i$ with $i<0$ by 1. Each replacement has the effect of reducing the number of variables in $g(\x)$ by one and dividing by $q^{k_+}$, so in the integral over the $B$-variables we must divide by $q^{k_+\,(k_--|B|)}$. Therefore to get the coefficients of the resulting integrals over large and small contours we must multiply the coefficient in (\ref{gint}) by $q^{-k_+\,(k_--|B|)}$. The result is
\[(-1)^{k_-+|B|}\,\t^{\s(B,\,S_-)-k_-\,|B|}\;{p^{|B|(|B|-1)/2}\ov q^{k_-(k_--1)/2+k_+\,(k_--|B|)}}\]
\[=(-1)^{k_-+|B|}\,\t^{-\s(S_-,\,B)}\;{p^{|B|(|B|+1)/2}\ov q^{k_-(k_--1)/2+k_+\,(k_--|B|)+|B|.}}\]
Here we used the identity
\be\s(U,V)+\s(V,U)=|U|\,|V|+|U\cap V|\label{UVVU}\ee
to obtain $\s(B,\,S_-)=-\s(S_-,\,B)+(k_-+1)\,|B|$.
Thus the integral on the right side of (\ref{PY3}) equals 
\[\sum_{B\subset S_-}(-1)^{k_-+|B|}\,\t^{-\s(S_-,\,B)}\;{p^{|B|(|B|+1)/2}\ov q^{k_-(k_--1)/2+k_+\,(k_--|B|)+|B|.}}\]
\be\times \int_{\C{r}^{|B|}}\int_{\C{R}^{k_+}}I(x,|B|+k_+,\x)\,\prod_{s_i\in B}\x_i^{-s_i}\,\prod_{i>0}\x_i^{-s_i}\,\prod_i d\x_i.\label{Bsum}\ee
When $|B|+k_-=0$ the integral is interpreted as 1. 

If we multiply the coefficient in (\ref{Bsum}) by the coefficient of the integral in (\ref{PY3}) and use
\[k(k-1)/2= k_-(k_--1)/2+k_+(k_+-1)/2+k_-k_+\]
and the bilinearity of $\s$ we get 
\[(-1)^{m+k_-+|B|}\t^{m(m-1)/2-km+\s(S_-,\,Y_-\bs B)+\s(S_+,\,Y)}\]
\[\times\,q^{k_+(k_+-1)/2+(k_+-1)\,|B|}\,p^{|B|\,(|B|+1)/2}\,\br{k-1}{m-1}.\]

Now instead of first summing over all $B\subset S_-$ and then over all $S_-\subset Y_-$, we reverse the order by first fixing $B$ and then summing over all $S_-$ satisfying $B\subset S_-\subset Y_-$. We change notation, so that in the end we get integrands like those in (\ref{PY3}) with $s_i\in S_-\cup S_+$. To do this we interchange $B$ and $S_-$, so the old $k_-$ becomes $|B|$ and the new $k_-$ is defined to be the new $|S_-|$. We label our indices so that $s_i\in S_-$ when $i<0$ and $s_i\in S_+$ when $i>0$. Then our formula becomes
\be\P_Y(x_m(t)\le x)=\sum_{k_\pm\ge0}\,\sum_{{S_\pm\subset Y_\pm\atop |S_\pm|=k_\pm}}c_{m,\,S_-,\,S_+}\int_{\C{R}^{k_+}}\int_{\C{r}^{k_-}} I(x,\,k,\,\x)\;\prod_{i} \x_i^{-s_i}\,d\x_i\label{P1}\ee 
where $k=k_-+k_+$ and 
\[c_{m,\,S_-,\,S_+}=(-1)^{m+k_-}\,\t^{m(m-1)/2-k_+m+\s(S_+,\,Y)}\,q^{k_+(k_+-1)/2+k_-\,(k_+-1)}\,p^{k_-\,(k_-+1)/2}\]
\[\times\,\sum_{Y_-\supset B\supset S_-}(-1)^{|B|}\,\t^{\s(B,\,Y_-\bs S_-)-m|B|}\,\br{|B|+k_+-1}{m-1}\]
\[=(-1)^{m+k_-}\,\t^{m(m-1)/2-k_+m+\s(S_+,\,Y)+k_-(k_-+1)/2}\,q^{k(k-1)/2}\]
\be\times\,\sum_{Y_-\supset B\supset S_-}(-1)^{|B|}\,\t^{\s(B,\,Y_-\bs S_-)-m|B|}\,\br{|B|+k_+-1}{m-1}.\label{CSS}\ee
When $k=0$ the integral in (\ref{P1}) is interpreted as 1. After our change of notation the restriction $k\ge1$ in the sum in (\ref{PY3}) becomes $|B|+k_+\ge1$.

To evaluate the sum in (\ref{CSS}) we use\footnote
{The left side equals $\sum_{1\le u_1<\cd<u_n\le|X|}\t^{u_1+\cd+u_n}$, which is the coefficient of $z^n$ in the expansion of $\prod_{k=1}^{|X|}(1+z\t^k)$ which, by the $\t$-binomial theorem \cite[p.26]{M}, equals $\sum_{n=0}^{|X|}\t^{n(n+1)/2}z^n\br{|X|}{n}$.}
\be\sum_{{U\subset X\atop|U|=n}}\t^{\s(U,\,X)}=\t^{n(n+1)/2}\,\br{|X|}{n}.\label{identity}\ee
We write the power of $\t$ in the sum in (\ref{CSS}) as
\[\s(B\bs S_-,\,Y_-\bs S_-)+\s(S_-,\,Y_-\bs S_-),\]
and then the sum in (\ref{CSS})  as
\[(-1)^{k_-}\t^{\s(S_-,\,Y_-\bs S_-)-mk_-}\,\sum_{Y_-\supset B\supset S_-}(-1)^{|B\bs S_-|}\,\t^{\s(B\bs S_-,\,Y_-\bs S_-)-m|B\bs S_-|}
\,\br{|B\bs S_-|+k-1}{m-1}.\]
The restriction $|B|+k_+\ge1$ is the same as $|B\bs S_-|+k\ge 1$.
If we apply (\ref{identity}) with $U=B\bs S_-$ and $X=Y_-\bs S_-$,
then the last sum may be written
\[\sum_{n=0}^{|Y_-\bs S_-|}\,(-1)^n\,\t^{n(n+1)/2-mn}\,\br{|Y_-\bs S_-|}{n}
\,\br{n+k-1}{m-1}\]
when $k>0$. When $k=0$ the sum begins with $n=1$. (Notice that all terms of the sum vanish unless $m\le|Y_-\bs S_-|+k$.)
The sum from $n=0$ is equal to\footnote{If we divide by the $n=0$ term of the series the result is the $\t$-hypergeometric function $_2\phi_1(a,b;c;\t,c/ab)$ with $a=\t^{-|Y_-\bs S_-|},\ b=\t^{k},\ c=\t^{k-m+1}$. An application of Heine's $q$-Gauss identity \cite[p. 14, eq. 1.52]{GR} with $q=\t$ gives the stated equality. We ``discovered'' the equality by a Maple computation. Doron Zeilberger \cite{Z} showed us that it follows from the $q$-Gauss identity.}
\[(-1)^{|Y_-\bs S_-|}\,\t^{|Y_-\bs S_-|\,(|Y_-\bs S_-|-2m+1)/2}\,
\br{k-1}{m-|Y_-\bs S_-|-1}.\]
When $k=0$ the condition $m\le|Y_-\bs S_-|+k$ is the same as $m\le|Y_-\bs S_-|$, and then the sum from $n=0$ equals zero. Therefore the sum from $n=1$ equals minus the $n=0$ term, which is
\be-\br{-1}{m-1}=(-1)^m\,\t^{-m(m-1)/2}.\label{k0}\ee

Using  
\[|Y_-\bs S_-|(|Y_-\bs S_-|+1)/2=\s(Y_-\bs S_-,\,Y_-\bs S_-)\]
 we see that when $k>0$ the sum in (\ref{CSS}) equals
\[(-1)^{|Y_-|}\,\t^{\s(Y_-,\,Y_-\bs S_-)-m\,|Y_-|}\,
\br{k-1}{m-|Y_-\bs S_-|-1}\]
and so 
\[c_{m,\,S_-,\,S_+}=(-1)^{m+|Y_-\bs S_-|}\,\t^{m(m-1)/2-k_+m+\s(S_+,\,Y)+\s(Y_-,\,Y_-\bs S_-)-m\,|Y_-|+k_-(k_-+1)/2}\]
\be\times\,q^{k(k-1)/2}\,\br{k-1}{m-|Y_-\bs S_-|-1}.\label{CmSS}\ee
When $k=0$ (i.e., $S_-=S_+=\emptyset$) the sum in (\ref{CSS}) equals $(-1)^m\,\t^{-m(m-1)/2}$ when $m\le|Y_-|$ and 0 otherwise, so $c_{m,\,\emptyset,\,\emptyset}=1$  when $m\le|Y_-|$ and 0 otherwise. Using (\ref{k0}) we find that this is exactly what we get if we substitute $k=0$ into (\ref{CmSS}). Therefore (\ref{CmSS}) holds for all $k_\pm\ge0$.

These are to be substituted into (\ref{P1})  to give $\P_Y(x_m(t)\le x)$. To obtain the formula for 
$\P_Y(x_m(t)=x)$ we replace $I(x,\,k,\,\x)$ in (\ref{P1}) by $\DI(x,\,k,\,\x)$, obtaining
\be\P_Y(x_m(t)=x)=\sum_{k_\pm\ge0}\,\sum_{{S_\pm\subset Y_\pm\atop |S_\pm|=k_\pm}}c_{m,\,S_-,\,S_+}\int_{\C{R}^{k_+}}\int_{\C{r}^{k_-}} \DI(x,\,k,\,\x)\;\prod_{i}\x_i^{-s_i}\,d\x_i.\label{P1+}\ee
Since the integral in (\ref{P1}) is interpreted as 1 for each $x$ when $k=0$, the integral here is interpreted as 0. Recall that in both formulas $s_i\in S_-$ when $i<0$ and $s_i\in S_+$ when $i>0$.

\begin{center}{\bf III. Formula for \mb$\P_Y(\e_x(t)=1)$}\end{center}

To obtain $\P_Y(\e_x(t)=1)$, the probability that site $x$ is occupied at time $t$ given the initial configuration $Y$, we replace the coefficients in (\ref{P1+}) by the sum of (\ref{CmSS}) over all $m$. The sum of those terms involving $m$ is
\[\sum_m(-1)^m\,\t^{m(m-1)/2-(k_++|Y_-|)\,m}\,\br{k-1}{m-|Y_-\bs S_-|-1},\]
and after making the substitution $m\to m+|Y_-\bs S_-|$ this becomes
\[(-1)^{|Y_-\bs S_-|}\,\t^{|Y_-\bs S_-|(|Y_-\bs S_-|-1)/2-(k_++|Y_-|)\,|Y_-\bs S_-|}\,\sum_m(-1)^m\,\t^{m(m-1)/2-mk}\,\br{k-1}{m-1}.\]
The factor outside the sum may be written
\[(-1)^{|Y_-\bs S_-|}\,\t^{-\s(Y_-\bs S_-,\,Y_-\bs S_-)-k\,|Y_-\bs S_-|},\]
while another application of the $\t$-binomial theorem shows that when $k\ge1$ the sum equals
\[(-1)^k\,\t^{-k(k+1)/2}\,\prod_{j=1}^{k-1}(1-\t^j).\]
If we use 
\be\s(Y_-,\,Y_-\bs S_-)-\s(Y_-\bs S_-,\,Y_-\bs S_-)=-\s(Y_-\bs S_-,\,S_-)+k_-\,|Y_-\bs S_-|,\label{sig1}\ee
\be\s(S_+,\,Y)=\s(S_+,\,Y_+\bs S_+)+k_+(k_++1)/2+k_+|Y_-\bs S_-|+k_+k_-,\label{sig2}\ee
which are obtained using the bilinearity of $\s$ and (\ref{UVVU}), and multiply by the factors in (\ref{CmSS}) not involving $m$, we obtain for the sum over $m$
\[c_{S_-,\,S_+}=(-1)^{k}\,\t^{\s(S_+,\,Y_+\bs S_+)-\s(Y_-\bs S_-,\,S_-)}\,q^{k(k-1)/2}\,\prod_{j=1}^{k-1}(1-\t^j)\]
\be=-\t^{\s(S_+,\,Y_+\bs S_+)-\s(Y_-\bs S_-,\,S_-)}\,\prod_{j=1}^{k-1}(p^j-q^j).\label{css}\ee
Thus, with these coefficients we have
\be\P_Y(\e_x(t)=1)=\sum_{k_\pm\ge0}\,\sum_{{S_\pm\subset Y_\pm\atop |S_\pm|=k_\pm}}
c_{\,S_-,\,S_+}\int_{\C{R}^{k_+}}\int_{\C{r}^{k_-}} \DI(x,\,k,\,\x)\;\prod_{i} \x_i^{-s_i}\,d\x_i,\label{PY}\ee
where the integral is interpreted as 0 when $k=0$.

\begin{center}{\bf IV. \mb$\P(\e_x(t)=1)$ for Bernoulli Initial Condition}\end{center}

We now derive a formula analogous to (\ref{PY}) for Bernoulli initial condition on $\Z$, with density $\r_-$ on $(\iy,\,0]$ and density $\r_+$ on $[1,\,\iy)$. It will be a sum over $k_\pm$ only. We shall begin with $\Z$ replaced by $(-M,\,N]$, apply (\ref{PY}), and then let $M,\,N\to\iy$. 

The probability for an initial configuration $Y=Y_-\cup Y_+$ is
\be\r_-^{|Y_-|}\,(1-\r_-)^{M-|Y_-|}\;\r_+^{|Y_+|}\,(1-\r_+)^{N-|Y_+|}.\label{initial}\ee
To obtain the coefficient of the multiple integral corresponding to the indices $k_\pm$ we have to multiply the above by $c_{S_-,\,S_+}$ as given by (\ref{css}), then sum over all $Y_\pm\supset S_\pm$, and finally multiply by $\prod \x_i^{-s_i}$ and sum over all $S_\pm$ with $|S_\pm|=k_\pm$. Because of the structure of $c_{S_-,\,S_+}$ we need only compute the two limits
\be\lim_{N\to\iy}\sum_{{S_+\subset [1,\,N]\atop |S_+|=k_+}}\sum_{{Y_+\supset S_+\atop Y_+\subset[1,\,N]}}\r_+^{|Y_+|}\,(1-\r_+)^{N-|Y_+|}\,\t^{\s(S_+,\,Y_+\bs S_+)}\,\prod_{i>0}\x_i^{-s_i},\label{limit1}\ee
where $S_+=\{s_i,\ld,s_{k_+}\}$ with $s_1<\cd<s_{k_+}$, and
\be\lim_{M\to\iy}\sum_{{S_-\subset (-M,\,0]\atop |S_-|=k_-}}\sum_{{Y_-\supset S_-\atop Y_+\subset(-M,\,0]}}\r_-^{|Y_-|}\,(1-\r_-)^{M-|Y_-|}\,\t^{-\s(Y_-\bs S_-,\,S_-)}\,\prod_{i<0}\x_i^{-s_i},\label{limit2}\ee
where $S_-=\{s_{-k_-},\ld,s_{-1}\}$ with $s_{-k_-}<\cd<s_{-1}$,
multiply them together, and then multiply by $\prod_{j=1}^{k-1}(p^j-q^j)$.
In (\ref{limit1}) the $|\x_i|$ are large and in (\ref{limit2}) the $|\x_i|$ are small.

The limit (\ref{limit1}) with the factor $\t^{\s(S_+,\,Y_+)}$ instead of $\t^{\s(S_+,\,Y_+\bs S_+)}$ was computed in \cite{TW3}.\footnote{The computation began with with formula (5) of the cited paper and the limit is the displayed expression preceding identity (8).} Using that result we find that (\ref{limit1}) is equal to
\be\prod_{i=1}^{k_+}{\r_+\ov \x_{i}\cd\x_{k_+}-1+\r_+(1-\t^{k_+-i+1})}.\label{ans1}\ee
We required that these $|\x_i|$ be sufficiently large to take the $N\to\iy$ limit.

For the limit (\ref{limit2}) we set
\[\tl s_i=-s_{-i}+1,\ \ \ \tl S_+=-S_-+1,\ \ \ \tl Y_+=-Y_-+1,\]
(observe that $\tl s_1<\cd<\tl s_{k_-}$) and use 
\[\s(Y_-\bs S_-,\,S_-)=\s(\tl S_+,\,\tl Y_+\bs\tl S_+).\]
If we set $\tl\x_i=\x_{-i}\inv,\ (i=1,\ld,k_-)$ then (\ref{limit2}) becomes
\[\lim_{M\to\iy}\sum_{{\tl S_+\subset [1,\ld,M]\atop |\tl S_+|=k_-}}\sum_{{\tl Y_+\supset \tl S_+\atop \tl Y_+\subset[1,\ld,M]}}\r_-^{|\tl Y_+|}\,(1-\r_-)^{M-|\tl Y_+|}\,\t^{-\s(\tl S_+,\,\tl Y_+\bs \tl S_+)}\,\prod_{i>0}\tl \x_i^{-\tl s_i+1}.\]
The limit formula (\ref{ans1}) with an obvious modification shows that (\ref{limit2}) equals
\[\prod_{i=1}^{k_-}{\r_-\,\tl\x_i\ov \tl\x_{i}\cd\tl\x_{k_-}-1+\r_-(1-\t^{i-k_--1})}=\prod_{i=1}^{k_-}{\r_-\,\x_{-i}\inv\ov (\x_{-i}\cd\x_{-k_-})\inv-1+\r_-(1-\t^{i-k_--1})}.\]
We require that these $|\x_i|$ be sufficiently small to take the $M\to\iy$ limit.

Denote by $\ph_\pm(k_\pm,\r_\pm,\x)$ the two limits, 
\[\ph_+(k_+,\r_+,\x)=\prod_{i=1}^{k_+}{\r_+\ov \x_{i}\cd\x_{k_+}-1+\r_+(1-\t^{k_+-i+1})},\]
\[\ph_-(k_-,\r_-,\x)=\prod_{i=1}^{k_-}{\r_-\,\x_{-i}\inv\ov (\x_{-i}\cd\x_{-k_-})\inv-1+\r_-(1-\t^{i-k_--1})}.\]
If we recall (\ref{css}) then we see that we have shown, formally, that that for Bernoulli on $\Z$ with densities~$\r_\pm$
\[\P(\e_t(x)=1)=-\sum_{k_\pm\ge0}\,\prod_{j=1}^{k-1}(p^j-q^j)\]
\be\times\int_{\C{R}^{k_+}}\int_{\C{r}^{k_-}} \DI (x,k,\,\x)\,\cdot\,\ph_-(k_-,\r_-,\x)\,\ph_+(k_+,\r_+,\x)\,\prod_{i}d\x_i\,.\label{siteprob1}\ee

The reason for the qualification ``formally'' is that we must justify the interchange of the limits $M,\,N\to\iy$ with the sum over $k_\pm$. To do this it is enough that we obtain a convergent series when the integrands in (\ref{css}) with their coefficients are replaced by their upper bounds over $M$ and $N$. Now not only was the limit (\ref{limit1}) determined in \cite{TW3}, but an expression was found for finite $N$, namely
\[\r_+^{k_+}\,\prod_{i=1}^{k_+}\x_i^{-i}\;\sum_{{t_i\ge0\atop\sum t_i\le N-k_+}}\prod_{i=1}^{k_+}\({1-\r_+(1-\t^{k-i+1})\ov \x_i\cd\x_{k_+}}\)^{t_i}.\] 
This gives the $N\to\iy$ limit if $R=|\x_i|$ is sufficiently large (depending on~$\t$) and also a uniform bound
\be A^{k_+^2}\,R^{-k_+^2/2+O(k_+)},\label{bound1}\ee
for some $A$ and arbitrarily large $R$. Similarly for (\ref{limit2}) we get a uniform bound
\be A^{k_-^2}\,r^{k_-^2/2-O(k_-)}\label{bound2}\ee
for some $A$ and arbitrarily small $r$.

As for the factor $\DI (x,k,\,\x)$ in the integrand, we have the easy bounds
\be\prod_{0<i<j}|f(\x_i,\,\x_j)|\le A^{k_+^2}\,R^{-k_+^2/2+O(k_+)},\label{bound3}\ee
\be\prod_{i<j<0}|f(\x_i,\,\x_j)|\le A^{k_-^2}\,r^{k_-^2/2-O(k_-)},\label{bound4}\ee
\be\prod_{{i<0\atop j>0}}|f(\x_i,\,\x_j)|\le A^{k_-k_+}\,r^{-k_-k_+}.\label{bound5}\ee
(For the last we used $rR\gg1$; see footnote \ref{rR}.) The rest of the integrand is at most $A^{k_+R+k_-r\inv}\,R^{O(k_+)}\,r^{-O(k_-)}$. If we use $k_-k_+\le k_-^2/2+k_+^2/2$ and combine the above bounds we get the bound
\[A^{k_+^2+k_-^2+k_+R+k_-r\inv}\,R^{-k_+^2+O(k_+)}\,r^{k_-^2/2-O(k_-)}.\]
If we take $R$ sufficiently large and $r$ sufficiently small (depending on this $A$) we get a convergent sum of integrals. So (\ref{siteprob1}) is justified.
\sp

\noindent{\bf Remark.} For SSEP the only nonzero terms in (\ref{siteprob1}) are those with $k_-=0,\,k_+=1$ and $k_+=0,\,k_-=1$, and we get a sum of two single integrals with $\r_\pm$ appearing only as coefficients:
\be\P(\e_t(x)=1)=\r_+\int_{\C{R}}{\x^{x-1}\,e^{\ep(\x)t}\ov\x-1}\,d\x+\r_-\int_{\C{r}}{\x^{x-1}\,e^{\ep(\x)t}\ov1-\x}\,d\x.\label{SSEP}\ee

\begin{center}{\it Symmetrization}\end{center}

In \cite{TW3}, where $\r_-=0$, we symmetrized the integrand in (\ref{siteprob1}), found that it was a determinant, and this led to a representation for $\P(x_m(t)\le x)$ as an integral of Fredholm determinants. Although it does not lead to as nice an expression, we can do an analogous symmetrization here.

Write the integrand in (\ref{siteprob1}) as
\be(1-\prod_i\x_i\inv)\,\prod_{{i<0\atop j>0}}f(\x_i,\,\x_j)\label{term1}\ee
\be\times \prod_{i<j<0}f(\x_i,\,\x_j)\;\prod_{i<0}{\x_i^{x}\,e^{\ep(\x_i)t}\ov 1-\x_i}\,\ph_-(k_-,\r_-,\x)\label{term2}\ee
\be\times \prod_{0<i<j}f(\x_i,\,\x_j)\;\prod_{i>0}{\x_i^{x}\,e^{\ep(\x_i)t}\ov 1-\x_i}\,\ph_+(k_+,\r_+,\x).\label{term3}\ee
 
The factor (\ref{term1}) is symmetric separately in the $\x_i$ with $i<0$ and the $\x_i$ with $i>0$. 

The factor (\ref{term3}) is a function of the $\x_i$ with $i>0$. Its symmetrization, given by formula (9) of \cite{TW3}, is
\[{1\ov k_+!}\,q^{k_+(k_+-1)/2}\,\prod_{i\ne j}f(\x_i,\,\x_j)\,\prod_i{\r_+\ov \x_i-1+\r_+(1-\t)}\;\prod_{i}{\x_i^{x}\,e^{\ep(\x_i)t}\ov 1-\x_i}.\]
(All indices positive.) By the identity (3) of \cite{TW2} this equals
\[{(-1)^{k_+}\ov k_+!}\,p^{-k_+(k_+-1)/2}\,\det\({1\ov p+q\x_i\x_j-\x_i}\)_{1\le i,\,j\le k_+}\,\prod_i{\r_+\,(q\x_i-p)\ov \x_i-1+\r_+(1-\t)}\;\x_i^{x}\,e^{\ep(\x_i)t}\]
\[={(-1)^{k_+}\ov k_+!}\,p^{-k_+(k_+-1)/2}\,\det(K_+(\x_i,\,\x_j))_{1\le i,\,j\le k_+},\]
where
\[K_+(\x,\x')={\r_+\,(q\x-p)\ov \x-1+\r_+(1-\t)}\;{\x^{x}\,e^{\ep(\x)t}\ov p+q\x \x'-\x}.\]

For the factors
\[\prod_{i<j<0}f(\x_i,\,\x_j)\;\ph_-(k_-,\r_-,\x),\]
we write as before $\tl \x_i=\x_{-i}\inv$, and this becomes
\[\prod_{i>j>0}{\tl\x_i-\tl\x_j\ov p \tl\x_i\tl\x_j+q-\tl\x_j}\;\prod_{i=1}^{k_-}{\r_-\,\tl\x_i\ov \tl\x_{i}\cd\tl\x_{k_+}-1+\r_-(1-\t^{i-k_--1})}\]
\[=\prod_{0<i<j}{\tl\x_j-\tl\x_i\ov q+ p \tl\x_i\tl\x_j-\tl\x_i}\;\prod_{i=1}^{k_-}{\r_-\,\tl\x_i\ov \tl\x_{i}\cd\tl\x_{k_+}-1+\r_-(1-\t^{i-k_--1})}\,.\]
The symmetrization of this is
\[{1\ov k_-!}p^{k_-(k_--1)/2}\,\;\prod_{i\ne j}{\tl\x_j-\tl\x_i\ov q+p\tl\x_i\tl\x_j- \tl\x_i}\,\prod_i{\r_-\,\tl\x_i\ov \tl\x_i-1+\r_-(1-\t\inv)}\ \ \ (\textrm{all indices positive})\]
\[={1\ov k_-!}\,p^{k_-(k_--1)/2}\,\;\prod_{i\ne j}{\x_j-\x_i\ov p+q \x_i \x_j-\x_i}\,\prod_i{\r_-\,\x_i\inv\ov \x_i\inv-1+\r_-(1-\t\inv)}\ \ \ (\textrm{all indices negative}).\]
Hence the symmetrization of (\ref{term2}) is (all indices negative)
\pagebreak
\[{1\ov k_-!}\,p^{k_-(k_--1)/2}\,\;\prod_{i\ne j}f(\x_i,\,\x_j)\,\prod_i{\r_-\,\x_i\inv\ov \x_i\inv-1+\r_-(1-\t\inv)}\,
\prod_{i}{\x_i^{x-1}\,e^{\ep(\x_i)t}\ov 1-\x_i}\]
\[={(-1)^{k_-}\ov k_-!}\,q^{-k_-(k_--1)/2}\,\det\({1\ov p+q\x_i\x_j-\x_i}\)\,\prod_i{\r_-\,(q-p\,\xi_i\inv)\ov \x_i\inv-1+\r_-(1-\t\inv)}\;\x_i^{x}\,e^{\ep(\x_i)t}\]
\[={(-1)^{k_-}\ov k_-!}\,q^{-k_-(k_--1)/2}\,\det(K_-(\x_i,\,\x_j))_{-k_-\le i,\,j\le -1},\]
where
\[K_-(\x,\x')={\r_-\,(q-p\,\x\inv)\ov \x\inv-1+\r_-(1-\t\inv)}\;{\x^{x}\,e^{\ep(\x)t}\ov p+q\x \x'-\x}.\]

So our formula may be written 
\[\P(\e_t(x)=1)=\sum_{{k_-,\,k_+\ge0\atop k_-+k_+\ge 1}}{1\ov k_-!\,k_+!}\,p^{-k_+(k_+-1)/2}\,q^{-k_-(k_--1)/2}\,\prod_{j=1}^{k_-+k_+-1}(q^j-p^j)\]
\[\times\int_{\C{R}^{k_+}}\int_{\C{r}^{k_-}}(1-\prod_i\x_i\inv)\,\prod_{{i<0\atop j>0}}f(\x_i,\,\x_j)\]
\be\times\,\det(K_-(\x_i,\,\x_j))_{-k_-\le i,\,j\le -1}\,\cdot\,\det(K_+(\x_i,\,\x_j))_{1\le i,\,j\le k_+}\,\cdot\,\prod_{i<0}d\x_i\,\cdot\,\prod_{i>0}d\x_i.\label{siteprob2}\ee

\noindent{\bf Remark 1}. Suppose the change of density occurs at $y$ rather than zero. Then to use the computations of (\ref{limit1}) and (\ref{limit2}) we would make the substitutions $Y_\pm=\tl Y_\pm+y,\ \ S_\pm=\tl S_\pm +y,\ \ s_i=\tl s_i+y$ and apply those computations to $\tl Y_\pm$, etc. The result in the end is that each $\x_i^x$ in the formula for $\DI(x,k,\,\x)$ becomes $\x_i^{x+y}$. In case $\r_-=\r_+$ it makes no difference what $y$ is so the formulas we got were independent of $x$, as they should be.

\noindent{\bf Remark 2}. To obtain formulas for the correlation function $\P(\e_t(x)=1,\,\e_0(0)=1)$ we write it as
\be\P(\e_t(x)=1)-\P(\e_t(x)=1,\,\e_0(0)=0).\label{difference}\ee
The first probability we know. For the second we modify (\ref{initial}) by multiplying by $1-\r_-$, the probability that site 0 is initially unoccupied, and take only those $Y_-$ that are contained in $(-\iy,-1]$. To use the preceding computation we use the one-one correspondence between subsets $S_-\subset(-\iy,-1]$ and subsets $S_-'\subset(-\iy,0]$ given by $S_-'=S_-+1$. Then $\s(Y_-\bs S_-,\,S_-)=\s(Y_-'\bs S_-',\,S_-')$ and, with obvious notation, $\prod \x_i^{-{s_i}}=\prod\x_i\times\prod \x_i^{-{s_i}'}$. It follows that for the second probability in (\ref{difference}) we multiply the integrands in (\ref{siteprob1}) and (\ref{siteprob2}) by $(1-\r_-)\,\prod_{i<0}\x_i$. Therefore for (\ref{difference}) itself we multiply the integrands by
\[1-(1-\r_-)\,\prod_{i<0}\x_i.\]

\begin{center}{\bf V. Formula for the Total Flux}\end{center}

Suppose $\max Y_-\le 0< \min Y_+$. When $t=0$ there are $|Y_-|$ particles $\le 0$, so $Q_t$, the total flux to the left across 0 at time $t$, is the number of particles $\le 0$ at time $t$ minus $|Y_-|$. Thus $\P_Y(Q_t\ge m)=\P_Y(x_{m+|Y_-|}(t)\le 0)$. Therefore in (\ref{P1}) and (\ref{CmSS}) we replace $m$ by $m+|Y_-|$ and set $x=0$.  If we use (\ref{sig1}) and (\ref{sig2}) and replace $m$ by $m+|Y_-|$ (so the new $m$ may be negative) then (\ref{CmSS}) becomes
\[c_{m+|Y_-|,\,S_-,\,S_+}=(-1)^{m+k_-}\,\t^{\s(S_+,\,Y_+\bs S_+)-\s(Y_-\bs S_-,\,S_-)+m(m-1)/2+k_+(k_++1)/2-mk_+}\]
\[\times\,q^{k(k-1)/2}\,\br{k-1}{m+k_--1}.\]
Using this, comparing with (\ref{css}), and following the argument of the last section we obtain the formulas 
\[\P(Q_t\ge m)=\sum_{k_\pm\ge0}(-1)^{m+k_-}\,\t^{m(m-1)/2+k_+(k_++1)/2-mk_+}\;q^{k(k-1)/2}\,\br{k-1}{m+k_--1}\]
\be\times\int_{\C{R}^{k_+}}\int_{\C{r}^{k_-}} I(0,k,\,\x)\,\cdot\,\ph_-(k_-,\r_-,\x)\,\ph_+(k_+,\r_+,\x)\,\prod_{i}d\x_i\,,\label{Tint}\ee
\[\P(Q_t\ge m)=\sum_{k_\pm\ge0}{(-1)^{m+k_+}\ov k_-!\,k_+!}\,\t^{m(m-1)/2-(m-1)k_+}\,q^{k_-k_+}\,\br{k-1}{m+k_--1}\]
\[\times\int_{\C{R}^{k_+}}\int_{\C{r}^{k_-}}\prod_{{i<0\atop j>0}}f(\x_i,\,\x_j)\,\cdot\,\det(K_-(\x_i,\,\x_j))_{-k_-\le i,\,j\le -1}\,\cdot\,\det(K_+(\x_i,\,\x_j))_{1\le i,\,j\le k_+}\,\cdot\,\prod_{i}d\x_i\,,\]
where in the expressions for $K_{\pm}(\x_i,\,\x_j)$ we set $x=0$.

For $\P(Q_t=m)$ we subtract from these what we get by replacing $m$ by $m+1$. To get a formula for $\langle e^{\la\,Q_t}\rangle$, the expected value of $e^{\la\, Q_t}$, we then multiply by $e^{\la\, m}$ and sum over all $m$. 

If $\P(Q_t\ge m)=p_m$ then $\langle e^{\la\, Q_t)}\rangle=\sum_{m=-\iy}^\iy (p_m-p_{m+1})\,e^{\la m}$. The series converges for all $\la$ and represents an entire function of $\la$.\footnote{For positive $m$, if $Q_t\ge m$ then the particle initially the $m$th to the right of 0 must have moved at least $m$ steps to the left at time $t$. The probability of this is less than the probability that a free particle would so move, and this probability is $O(e^{-m\,\log m+O(m)})$ as $m\to+\iy$. Similarly $\P(Q_t\ge~m)=1-O(e^{-|m|\,\log| m|+O(|m|)})$ as $m\to-\iy$.} Therefore if we find a formula that holds when $\Re\,\la>0$, and the formula represents an entire function of $\la$, it will hold in general. When $\Re\,\la>0$ we can write the series as the difference of two convergent series in the obvious way. Changing the summation variable in the second series and combining gives $(1-e^{-\la})\,\sum_{m=-\iy}^\iy p_m\,e^{\la m}$. This is what we do below. 

The factor that involves $m$ in the first formula is
\[(-1)^{m}\,\t^{m(m-1)/2- m\,k_+}\,\br{k-1}{m+k_--1}\]
When we make the replacement, subtract, multiply by $e^{\la\, m}$, and sum over $m$ we get
\[(1-e^{-\la})\,\sum_{m=-k_-+1}^{k_+}\,(-1)^{m}\,\t^{m(m-1)/2-m\,k_+}\,e^{\la\,m}\,\br{k-1}{m+k_--1}\]
and another application of the $\t$-binomial theorem shows that this equals
\[(-1)^{k_+}\,(1-e^{-\la})\,e^{\la\,k_+}\,\t^{-k_+(k_++1)/2}\,\prod_{j=1}^{k-1}(1-e^{-\la}\,\t^j)\]
\[=(-1)^{k_+}\,e^{\la\,k_+}\,\t^{-k_+(k_++1)/2}\,\prod_{j=0}^{k-1}(1-e^{-\la}\,\t^j).\]
Therefore
\[\langle e^{\la\, Q_t)}\rangle=\sum_{k_\pm\ge0}(-1)^{k}\,q^{k(k-1)/2)}\,e^{\la\,k_+}\,\prod_{j=0}^{k-1}(1-e^{-\la}\,\t^j)\]
\[\times\int_{\C{R}^{k_+}}\int_{\C{r}^{k_-}} I(0,k,\,\x)\,\cdot\,\ph_-(k_-,\r_-,\x)\,\ph_+(k_+,\r_+,\x)\,\prod_{i}d\x_i.\]

To justify the interchange of the sums over $m$ and $k_\pm$, and to show that the sum represents an entire function of $\la$, we use estimates analogous to those used in Sec.~IV. For the integral (\ref{Tint}) we have the bound $A^{k^2+k_+R+k_-r\inv}\,R^{-k_+^2+O(k_+)}\,r^{k_-^2/2-O(k_-)}$ and for its coefficient another bound $A^{k^2}$ for all $m$, since $|m|\le k$ when the coefficient is nonzero. And $|e^{\la m}|\le e^{k\,\Re\la}$ for such $m$. If $R$ is large enough and $r$ small enough the sum (over $k_\pm$ and $m$) of the bounds is finite and uniformly bounded in bounded $\la$-sets. This justifies the computation.

In the case of SSEP ($\t=1$) our formula simplifies since
\[\ph_\pm(k_\pm,\r_\pm,\x)=\r_\pm^{k_\pm}\,\ph_\pm(k_\pm,1,\x),\]
and we can write it as
\pagebreak
\[\langle e^{\la\, Q_t}\rangle=\sum_{k_\pm\ge0}(-1)^{k}\,2^{-k(k-1)/2}\,\r_-^{k_-}\,\r_+^{k_+}\,(1-e^{-\la})^{k_-}\,
(e^{\la}-1)^{k_+}\]
\[\times\int_{\C{R}^{k_+}}\int_{\C{r}^{k_-}} I(0,k,\,\x)\,\cdot\,\ph_-(k_-,1,\x)\,\ph_+(k_+,1,\x)\,\prod_{i}d\x_i\,.\]
In particular, this is a function of the two quantities $\r_-\,(1-e^{-\la})$ and $\r_+\,(e^{\la}-1)$.

This fact was derived in \cite{DG} in a different way. The authors went further to show that it is a function of a particular combination of these quantities. This allowed them to reduce to the case when $\r_-=0$ and to obtain, for general $\r_\pm$, a series  for the generating function involving integrals only over small contours. 

\begin{center}{\bf Acknowledgments}\end{center}

We thank Doron Zeilberger for help with the combinatorial identity in Sec.~III. 

This work was supported by the National Science Foundation through grants DMS-0906387 (first author) and DMS-0854934 (second author).

\end{document}